\documentclass[letterpaper, 10 pt, conference]{ieeeconf}  

\usepackage{graphicx} 
\usepackage{amsmath,bm,times} 
\usepackage{amssymb}  
\usepackage{bm}
\usepackage{tikz}
\usetikzlibrary{shapes,arrows,backgrounds,fit,positioning}
\usepackage{subfigure}
\usepackage{cite}
\usepackage{booktabs} 

\newtheorem{mydef}{\textbf{Definition}} [section]
\newtheorem{assump}[mydef] {\textbf{Assumption}}
\newtheorem{thm}[mydef] {\textbf{Theorem}}

\newtheorem{prop}[mydef]{\textbf{Proposition}}
\newtheorem{cor}[mydef]{\textbf{Corollary}}

\newtheorem{ex}[mydef]{\textbf{Example}}

\renewcommand{\themydef}{\thesection.\arabic{mydef}}

\DeclareMathOperator*{\argmax}{arg\,max}

\title{\LARGE \bf 
	Passivity and Evolutionary Game Dynamics}

\author{Shinkyu Park, Jeff S. Shamma, and Nuno C. Martins
  \thanks{The work of S. Park and N. C. Martins was supported by AFOSR grant FA95501510367, NSF ECCS 1408320 and NSF: ECCS 1135726. The work of J. S. Shamma was supported by funding from King Abdullah University of Science and Technology (KAUST).}
  \thanks{Shinkyu Park is with Senseable City Laboratory and Computer Science and Artificial Intelligence Laboratory, Massachusetts Institute of Technology, Cambridge, MA 02139, USA. {\tt shinkyu@mit.edu}}
  \thanks{Nuno C. Martins is with the Department of Electrical and Computer Engineering, University of Maryland College Park, College Park, MD 20742-4450, USA. {\tt nmartins@umd.edu}}
  \thanks{Jeff S. Shamma is with King Abdullah University of Science and Technology (KAUST), Computer, Electrical and Mathematical Science and Engineering Division (CEMSE), Thuwal 23955--6900, Saudi Arabia. {\tt jeff.shamma@kaust.edu.sa}}
}

\IEEEoverridecommandlockouts
\begin{document}

\maketitle

\begin{abstract}
  This paper investigates an energy conservation and dissipation -- passivity -- aspect of dynamic models in evolutionary game theory. We define a notion of passivity using the state-space representation of the models, and we devise systematic methods to examine passivity and to identify properties of passive dynamic models. Based on the methods, we describe how passivity is connected to stability in population games and illustrate stability of passive dynamic models using numerical simulations.
\end{abstract}

\section{Introduction} \label{sec_introduction}
Of central interest in evolutionary game theory \cite{weibull1995_mit, hofbauer2003_ams} is the study of strategic interactions among players in large populations. Each player engaged in a game chooses a strategy among a finite set of options and repeatedly revises its strategy choice in response to given payoffs. The distribution of strategy choices by the players define the population state, and evolutionary dynamics represent the strategy revision process and prompt the time-evolution of the population state.

One of main focuses of the study is on analyzing the population state trajectories induced by evolutionary dynamics to identify asymptotes of the trajectories and establishing \textit{stability} in population games. This paper contributes to the study by investigating a passivity -- abstraction of energy conservation and dissipation -- aspect of evolutionary dynamic models (EDMs) specifying the dynamics and by demonstrating how passivity yields stability in population games.

The most relevant work in literature investigates various aspects of evolutionary dynamics and associated stability concepts in population games. To mention a few, Brown and von Neumann \cite{brown1950_ams} studied \textit{Brown-von Neumann-Nash (BNN)} dynamics to examine the existence of optimal strategy choices in a zero-sum two-player game. Taylor and Jonker \cite{taylor1978_mb} studied \textit{replicator} dynamics and established a connection between \textit{evolutionarily stable strategies} \cite{maynard_smith1973} and stable equilibrium points of replicator dynamics. Later the result was strengthened by Zeeman \cite{zeeman1980_gtds} who also proposed a stability concept for games under replicator dynamics. Gilboa and Matshu \cite{gilboa1991_econometrica} considered \textit{cyclic stability} in games under best-response dynamics.

In succeeding work, stability results are established using broader classes of evolutionary dynamics. Swinkels \cite{swinkels1993_geb} considered a class of \textit{myopic adjustment} dynamics and studied \textit{strategic stability} of Nash equilibria under these dynamics. Ritzberger and Weibull \cite{ritzberger1995_econometrica} considered a class of \textit{sign-preserving selection} dynamics and studied asymptotic stability of faces of the population state space under these dynamics.

In a recent development of evolutionary game theory, Hofbauer and Sandholm \cite{hofbauer2009_jet} proposed \textit{stable games} and established stability of Nash equilibria of these games under a large class of evolutionary dynamics. The class includes \textit{excess payoff/target (EPT)} dynamics, \textit{pairwise comparison} dynamics, and \textit{perturbed best response (PBR)} dynamics. Fox and Shamma \cite{fox2013_games} later revealed that the aforementioned class of evolutionary dynamics exhibit passivity. Based on passivity methods from dynamical system theory, the authors established $\mathbb L_2$-stability of evolutionary dynamics in a generalized class of stable games. Mabrok and Shamma \cite{7799211} discussed passivity for higher-order evolutionary dynamics. In particular, the authors investigated a connection between passivity of linearized dynamics and stability in generalized stable games using robust control methods. Passivity methods are also adopted in games over networks to analyze stability of Nash equilibrium \cite{gadjov2017_arxiv}.

Inspired by passivity analysis presented in \cite{fox2013_games}, we investigate the notion of passivity for EDMs of evolutionary dynamics more in-depth. \textbf{Our main goals are} \textbf{(i)} to define passivity of EDMs that admit realizations in a finite-dimensional state space; \textbf{(ii)} to develop methods to examine passivity and to identify various properties of passive EDMs; and \textbf{(iii)} based on the above results, to establish a connection between passivity and stability in population games.

\subsection{Summary of the Main Contributions}
\begin{enumerate}
\item We define and characterize the notions of $\delta$-passivity for EDMs that admit state-space representation. Based on the characterization result, we demonstrate how to examine $\delta$-passivity of EDMs for the replicator dynamics, EPT dynamics, pairwise comparison dynamics, and PBR dynamics.

\item We investigate certain properties of $\delta$-passive EDMs with respect to payoff monotonicity and total payoff perturbation.
  
\item Based on the above results, we establish a connection between $\delta$-passivity and asymptotic stability of Nash equilibria in population games. We also illustrate the stability using numerical simulations.
\end{enumerate}

\subsection{Paper Organization}
In Section \ref{sec:background}, we present background materials on evolutionary game theory that are needed throughout the paper. In Section \ref{sec:passive_evo_dyn}, we define $\delta$-passivity of EDMs and present an algebraic characterization of $\delta$-passivity in terms of the vector field which describes state-space representation of EDMs. In Section \ref{sec:properties}, using the characterization result, we investigate properties of $\delta$-passive EDMs. In Section \ref{sec:numerical_ex}, we establish a connection between $\delta$-passivity and asymptotic stability in population games, and present some numerical simulation results to illustrate the stability. We end the paper with conclusions in Section \ref{sec:conclusion}.

\subsection{Notation}
\begin{itemize}
\item $a_i$ -- given a vector $a$ in $\mathbb R^n$, we denote its $i$-th entry as $a_i$.
  
\item $[a]_+$, $[a_i]_+$ -- the non-negative part of a vector ${a = \begin{pmatrix} a_1 & \cdots & a_n \end{pmatrix}^T}$ in $\mathbb R^n$ and its $i$-th entry defined by $[a]_+ \overset{\mathrm{def}}{=} \begin{pmatrix} [a_1]_+ & \cdots & [a_n]_+\end{pmatrix}^T$ and ${[a_i]_+ \overset{\mathrm{def}}{=} \max\{a_i, 0\}}$, respectively.

\item $\mathbf 1, I, e_i$ -- the vector with all entries equal to $1$, the identity matrix, and the $i$-th column of $I$, respectively.\footnote{We omit the dimensions of vectors and matrices whenever they are clear from context.}

\item $\nabla_x f, \nabla_x^2 f$ -- the gradient and Hessian of a real-valued function $x \mapsto f(x)$ with respect to $x$, respectively, provided that they exist.

\item $D \mathcal F$ -- the differential of a mapping $\mathcal F: \mathbb X \to \mathbb R^n$.
  
\item $\mathrm{int}(\mathbb X ), \mathrm{bd}(\mathbb X)$ -- the (relative) interior and boundary of a set $\mathbb X$ in its affine hull, respectively.

\item $\mathbb R_+^n \left( \mathbb R_-^n \right)$ -- the set of $n$-dimensional element-wise non-negative (non-positive) vectors. For $n=1$, we omit the superscript $n$ and adopt $\mathbb R_+ \left( \mathbb R_- \right)$.
  
\item $\| \cdot \|$ -- the Euclidean norm.

\end{itemize}

\section{Background on Evolutionary Game Theory} \label{sec:background}
Consider a population of players engaged in a game where each player selects a (pure) strategy from the set of available strategies denoted by $\{1, \cdots, n\}$.\footnote{Population games, in general, account multiple populations of players, and the strategy sets are allowed to be distinct across populations (see \cite{hofbauer2009_jet} for details). However, for simple and clear presentation, we restrict our attention to a single-population case.} Suppose that the population consists of a continuum of players. The population states, which describe the distribution of strategy choices by players, constitute a simplex defined by ${\mathbb X \overset{\mathrm{def}}{=} \{x \in \mathbb R_+^n \,|\, \sum_{i=1}^n x_i = 1 \}}$ where $x_i$ is a fraction of population choosing strategy $i$. Let us denote the tangent space of $\mathbb X$ by $\mathbb {TX} = \left\{ z \in \mathbb R^n \,|\, \sum_{i=1}^n z_i = 0 \right\}$.

In population games, an element $p$ in $\mathbb R^n$, so-called the payoff vector, is assigned to the population state $x$ by payoff functions where the $i$-th entry $p_i$ of $p$ represents the payoff given to the players choosing strategy $i$. In this work, we mainly focus on investigating passivity of evolutionary dynamic models and briefly demonstrate how passivity yields stability in population games. For the later purpose, we consider the payoff functions given below:
\begin{enumerate}
\item \textit{Static Payoff \cite{hofbauer2009_jet}:}
  \begin{align} \label{eq:static_payoff}
    p(t) = \mathcal F(x(t))
  \end{align}
  
\item \textit{Smoothed Payoff \cite{fox2013_games}:}
  \begin{align} \label{eq:smoothed_payoff}
    \dot p(t) = \lambda (\mathcal F(x(t)) - p(t))
  \end{align}
\end{enumerate}

In \eqref{eq:smoothed_payoff}, we note that at the equilibrium points, it holds that $p(t) = \mathcal F(x(t))$. Throughout the paper, we adopt the definition of a Nash equilibrium as follows.
\begin{mydef}
  Given a payoff function $\mathcal F: \mathbb X \to \mathbb R^n$, the population state $x^{\mathrm{NE}} \in \mathbb X$ is a Nash equilibrium of $\mathcal F$ if it holds that
  \begin{align*}
    x_i^{\mathrm{NE}}>0 \text{ implies } i \in \argmax_{j \in \{1, \cdots, n\}} \mathcal F_j(x^{\mathrm{NE}})
  \end{align*}
\end{mydef}

\subsection{Evolutionary Dynamic Models} \label{sec:evo_dyn}
Evolutionary dynamics describe how the population state evolves over time in response to given payoffs and are specified by evolutionary dynamic models (EDMs). We focus on EDMs  that admit state-space representation of the following form:
\begin{align} \label{eq:edm}
  \dot x(t) = \mathcal V(p(t), x(t)), ~ x(0) \in \mathbb X
\end{align}
where $p(t)$, $x(t)$, and $\dot x(t)$ take values in $\mathbb R^n$, $\mathbb X$, and $\mathbb {TX}$, respectively. We assume that the vector field $\mathcal V: \mathbb R^n \times \mathbb X \to \mathbb {TX}$ is \textit{well-defined} in a sense that it continuously depends on $(p(t), x(t))$, and for each initial value $x(0)$ in $\mathbb X$ and each payoff vector trajectory $p(\cdot): \mathbb R_+ \to \mathbb R^n$ in $\mathfrak P$, there exists a unique solution $x(\cdot): \mathbb R_+ \to \mathbb X$ to \eqref{eq:edm} that belongs to $\mathfrak X$, where $\mathfrak P$ and $\mathfrak X$ are defined by

{\small
  \begin{align*}
    \mathfrak P &\overset{\mathrm{def}}{=} \left\{ p (\cdot) \,\bigg|\, p(t) \in \mathbb X \text{ and } \int_0^t \|\dot p(\tau)\|^2 \, \mathrm d \tau < \infty \text{ for $t \in \mathbb R_+$} \right\} \\
    \mathfrak X &\overset{\mathrm{def}}{=} \left\{ x (\cdot) \,\bigg|\, x(t) \in \mathbb X \text{ and } \int_0^t \|\dot x(\tau)\|^2 \, \mathrm d \tau < \infty \text{ for $t \in \mathbb R_+$} \right\}
  \end{align*}}

Throughout the paper, we assume that EDMs satisfy the following path-connectedness condition.
\begin{assump} [Path-Connectedness] \label{assump:path_connected}
  We define the set of equilibrium points of \eqref{eq:edm} by ${\mathbb S \overset{\mathrm{def}}{=} \{ (p, x) \in \mathbb R^n \times \mathbb X \,|\, \mathcal V(p,x) = \mathbf 0 \}}$ and its projection on $\mathbb R^n \times \{x\}$ by $\mathbb S_x \overset{\mathrm{def}}{=} \{ p \in \mathbb R^n \,|\, (p, x) \in \mathbb S \}$ for each $x$ in $\mathbb X$. We assume that the set $\mathbb S_x$ is path-connected for every $x$ in $\mathbb X$. In other words, for every $p_0, p_1$ in $\mathbb S_x$, there exists a piece-wise smooth path from $p_0$ to $p_1$ that is contained in $\mathbb S_x$. \hfill $\square$
\end{assump}

We list below a few examples of EDMs found in the literature.
\begin{ex} \label{ex:evo_dynamics}
  For each $i$ in $\{1, \cdots, n\}$,
  \begin{enumerate}
  \item \textit{Replicator Dynamics \cite{taylor1978_mb}: } 
    \begin{align}
      \dot{x}_i(t) = x_i(t) \left( p_i(t) - p^T(t)x(t)\right) \label{eq:rep}
    \end{align}

  \item \textit{BNN Dynamics \cite{brown1950_ams}: } 
    \begin{align}
      \dot{x}_i(t) = [\hat{p}_i(t)]_{+} - x_i(t) \sum_{j=1}^{n} [\hat{p}_j(t)]_{+} \label{eq:bnn}
    \end{align}

   \item \textit{Smith Dynamics \cite{doi:10.1287/trsc.18.3.245}: } 
     \begin{align}
       \dot{x}_i(t) &= \sum_{j=1}^{n} x_j(t) [p_i(t) - p_j(t)]_{+} \nonumber \\
                    &\quad- x_i(t) \sum_{j=1}^{n} [p_j(t) - p_i(t)]_{+} \label{eq:smith}
     \end{align}

   \item \textit{Logit Dynamics \cite{hofbauer2007_jet}: } 
     \begin{align}
       \dot x_i(t) = \frac{\exp(\eta^{-1} p_i(t))}{\sum_{j=1}^{n} \exp(\eta^{-1} p_j(t))} - x_i(t) \label{eq:logit}
     \end{align}
  \end{enumerate}
  where $\hat p = \begin{pmatrix} \hat p_1 & \cdots & \hat p_n \end{pmatrix}^T$ in \eqref{eq:bnn} is the \textit{excess payoff vector} defined as $\hat{p} = p - p^T x \mathbf{1}$ and $\eta$ in \eqref{eq:logit} is a positive constant. \hfill $\square$
\end{ex}

Associated with Assumption \ref{assump:path_connected}, it can be verified that all the EDMs described in Example \ref{ex:evo_dynamics} satisfy the path-connectedness assumption. For instance, consider the set $\mathbb S = \{(p,x) \,|\, x \in \operatorname*{arg\,max}_{y \in \mathbb{X}} p^Ty\}$ which is the set of equilibrium points of \eqref{eq:bnn} and \eqref{eq:smith}. It can be verified that if both $(p_0,x)$ and $(p_1,x)$ belong to $\mathbb S$ then so does ${(\lambda p_0 + (1-\lambda) p_1, x)}$ for all $\lambda$ in $[0,1]$.

\subsection{Passivity and Stability}
As it will be presented in Section \ref{sec:passive_evo_dyn}, passivity of EDMs allows us to construct an energy function which can be used to identify asymptotes of population state trajectories induced by EDMs in a class of population games. This observation naturally explains asymptotic stability of equilibrium points of EDMs.

After we present our main results on passivity, we show how passivity is connected to stability in population games and provide numerical simulation results to illustrate the stability. We defer to a future extension of this paper for in-depth stability analysis for passive EDMs in a large class of population games.

\section{$\delta$-Passivity of Evolutionary Dynamic Model} \label{sec:passive_evo_dyn}
We proceed with defining $\delta$-passivity of EDMs, and then characterize $\delta$-passivity conditions in terms of the vector field $\mathcal V$ in \eqref{eq:edm}. Using the characterization, we examine $\delta$-passivity of representative EDMs in the literature and explore important properties of $\delta$-passive EDMs.

\subsection{Definition of $\delta$-Passivity}
Given an EDM \eqref{eq:edm} with a payoff vector trajectory $p(\cdot) \in \mathfrak P$ and a resulting population state trajectory $x(\cdot) \in \mathfrak X$, consider the following inequality: For $t \geq t_0 \geq 0$,
\begin{align} \label{eq:passivity_inequality}
  &\int_{t_0}^t \left[ \dot p^T(\tau) \dot x(\tau) - \eta \dot x^T(\tau) \dot x (\tau) \right] \,\mathrm d \tau \nonumber \\
  &\geq \mathcal S(p(t), x(t)) - \mathcal S(p(t_0), x(t_0))
\end{align}
where $\mathcal S: \mathbb R^n \times \mathbb X \to \mathbb R_+$ is a $\mathcal C^1$ function  and $\eta$ is a non-negative constant. Using \eqref{eq:passivity_inequality}, we state the definition of $\delta$-passivity as follows.
\begin{mydef} \label{def:passivity}
  Given an EDM \eqref{eq:edm}, we consider the following two cases:
  \begin{enumerate}
  \item  The EDM is \textit{$\delta$-passive} if there is a $\mathcal C^1$ function $\mathcal S$ for which \eqref{eq:passivity_inequality} holds with $\eta = 0$ for every $p(\cdot)$ in $\mathfrak P$.
  \item The EDM is \textit{strictly output $\delta$-passive} with index $\eta>0$ if there is a $\mathcal C^1$ function $\mathcal S$ for which \eqref{eq:passivity_inequality} holds for every $p(\cdot)$ in $\mathfrak P$.\footnote{It immediately follows from Definition \ref{def:passivity} that strict output $\delta$-passivity implies $\delta$-passivity.}
  \end{enumerate}  
\end{mydef}
We refer to $\mathcal S$ as a \textit{storage function}. The function $\mathcal S$ describes the energy stored in the EDM, and the $\delta$-passivity inequality \eqref{eq:passivity_inequality} suggests that the variation in the stored energy $\mathcal S(p(t), x(t)) - \mathcal S(p(t_0), x(t_0))$ is upper bounded by the supplied energy $\int_{t_0}^t [ \dot p^T(\tau) \dot x(\tau) - \eta \dot x^T(\tau) \dot x(\tau) ] \,\mathrm d\tau$.

We adopt the following definition of a strict storage function.
\begin{mydef}
  Let $\mathcal S: \mathbb R^n \times \mathbb X \to \mathbb R_+$ be a storage function of a $\delta$-passive EDM \eqref{eq:edm}. We call $\mathcal S$ \textit{strict} if the following condition holds: For all $(p,x)$ in $\mathbb R^n \times \mathbb X$,
  \begin{align}
    \nabla_x^T \mathcal S(p, x) \mathcal V(p,x) = 0 \text{ if and only if } \mathcal V(p,x)=\mathbf 0
  \end{align}
\end{mydef}

\subsection{Characterization of $\delta$-Passivity Condition} \label{sec:algebraic_conditions}
Let us consider the following relations: For $(p,x)$ in ${\mathbb R^n \times \mathbb X}$,
\begin{align}
  \nabla_p \mathcal S(p, x) &= \mathcal V(p, x) \tag{\textbf{P1}} \label{eq:p_condition_01} \\
  \nabla_x^T \mathcal S(p, x) \mathcal V(p, x) & \leq - \eta \mathcal V^T(p, x) \mathcal V(p, x)  \tag{\textbf{P2}} \label{eq:p_condition_02}
\end{align}
where $\mathcal S: \mathbb R^n \times \mathbb X \to \mathbb R_+$ is a $\mathcal C^1$ function, $\mathcal V: \mathbb R^n \times \mathbb X \to \mathbb T\mathbb X$ is the vector field in \eqref{eq:edm}, and $\eta$ is a non-negative constant. In the following theorem, using \eqref{eq:p_condition_01} and \eqref{eq:p_condition_02}, we characterize the $\delta$-passivity condition.

\begin{thm} \label{thm:passivity_edm}
  For a given EDM \eqref{eq:edm}, the following two statements are true:
  \begin{enumerate}
  \item  The EDM is \textit{$\delta$-passive} if there is a $\mathcal C^1$ function $\mathcal S$ satisfying \eqref{eq:p_condition_01} and \eqref{eq:p_condition_02} with $\eta = 0$.
  \item The EDM is \textit{strictly output $\delta$-passive} with index $\eta>0$ if there is a $\mathcal C^1$ function $\mathcal S$ satisfying \eqref{eq:p_condition_01} and \eqref{eq:p_condition_02}.
  \end{enumerate}
\end{thm}

The proof of Theorem \ref{thm:passivity_edm} is given in Appendix \ref{proof:passivity_edm}.

The following are interpretations of the conditions \eqref{eq:p_condition_01} and \eqref{eq:p_condition_02}. The condition \eqref{eq:p_condition_01} ensures that the integral $\int_P \mathcal V(\mathbf p, x) \cdot \mathrm d\mathbf p$ does not depend on the choice of the path $P$. This condition is an important requirement in establishing stability \cite{sandholm2013_dga}. The condition \eqref{eq:p_condition_02} implies that with the payoff vector fixed, the population state $x(t)$ evolves along a trajectory for which the function $\mathcal S$ decreases. This condition plays an important role in identifying \text{stable} equilibria of EDMs.

\subsection{Assessment of $\delta$-Passivity}
Using Theorem \ref{thm:passivity_edm}, we evaluate $\delta$-passivity of the following four EDMs. The definition of a revision protocol given below is used to describe some of the EDMs.

\begin{mydef} \label{def:revision_proc}
  A function $\varrho = ( \varrho_1 ~ \cdots ~ \varrho_n )^T$ is called the \textit{revision protocol} in which each entry $\varrho_i$ is defined as $\varrho_i: \mathbb R^n \to \mathbb R_+$.
\end{mydef}

\subsubsection{Replicator Dynamics \cite{taylor1978_mb}} For each $i$ in $\{1, \cdots, n\}$,
\begin{align} \label{eq:rep_dyn} 
  \dot x_i(t) = x_i(t) \left( p_i(t) - p^T(t) x(t) \right)
\end{align}

\begin{prop} \label{prop:rep_dyn}
  The EDM \eqref{eq:rep_dyn}  is not $\delta$-passive.
\end{prop}

\subsubsection{Excess Payoff/Target (EPT) Dynamics \cite{sandholm2005_jet}} For each $i$ in $\{1, \cdots, n\}$, 
\begin{align}
  \dot x_i(t) = \varrho_i (\hat p(t)) - x_i(t) \mathbf 1^T \varrho(\hat p(t)) \label{eq:ept_dyn}
\end{align}
where $\hat p$ is the excess payoff vector defined as $\hat p = p - p^T x \mathbf 1$, and the revision protocol $\varrho = ( \varrho_1 ~ \cdots ~ \varrho_n )^T$ satisfies the following two conditions -- Integrability \eqref{eq:integrability} and Acuteness \eqref{eq:acuteness}:
\begin{align}
  & \nabla_{\hat p} \gamma(\hat p) =  \varrho(\hat p) \tag{\textbf I} \label{eq:integrability} \\
  & \hat p^T \varrho (\hat p) > 0 \text{ if } \hat p \in \mathbb R^n \setminus \mathbb R^n_- \tag{\textbf A} \label{eq:acuteness}
\end{align}
where $\gamma: \mathbb R^n \to \mathbb R$ is a $\mathcal C^1$ function. Note that \eqref{eq:bnn}  is a particular case of \eqref{eq:ept_dyn}. The following proposition establishes $\delta$-passivity of \eqref{eq:ept_dyn}.

\begin{prop} \label{prop:ept_passivity}
  The EDM \eqref{eq:ept_dyn}  is $\delta$-passive and has a strict storage function given by $\mathcal S(p,x) = \gamma(\hat p)$ where $\gamma: \mathbb R^n \to \mathbb R_+$ is a non-negative $\mathcal C^1$ function satisfying \eqref{eq:integrability}.\footnote{The condition \eqref{eq:integrability} only ensures the existence of a $\mathcal C^1$ function $\gamma$ which would be negative. However, in the proof of Proposition \ref{prop:ept_passivity}, based on Assumption \ref{assump:path_connected} we show that there is a choice of $\gamma$ that is non-negative and satisfies \eqref{eq:integrability}.}
\end{prop}

\subsubsection{Pairwise Comparison Dynamics \cite{sandholm2010_games}} For each $i$ in $\{1, \cdots, n\}$, 
\begin{align}
  \dot x_i(t) &= \sum_{j=1}^nx_j(t) \varrho_i (p_i(t) - p_j(t)) \nonumber \\
              &\quad - x_i(t) \sum_{j=1}^n \varrho_j (p_j(t) - p_i(t)) \label{eq:pc_dyn}
\end{align}
The revision protocol $\varrho = (\varrho_1 ~ \cdots ~ \varrho_n )^T$ satisfies the following condition -- Sign Preservation \eqref{eq:sign_preservation}:\footnote{The function $\mathrm{sgn}: \mathbb{R} \to \{-1, 0, 1\}$ is defined as ${\mathrm{sgn}(a) \overset{\mathrm{def}}{=} \begin{cases} 1 & \text{ if } a>0 \\ 0 & \text{ if } a=0 \\ -1 & \text{ if } a<0 \end{cases}}$.}
\begin{align}
  \mathrm{sgn} ( \varrho_{i}(p_i - p_j) ) = \mathrm{sgn} ( [ p_i - p_j ]_{+} ) \tag{\textbf{SP}} \label{eq:sign_preservation}
\end{align}
Note that \eqref{eq:smith} is a particular case of \eqref{eq:pc_dyn}. The following proposition establishes $\delta$-passivity of \eqref{eq:pc_dyn}
\begin{prop} \cite{fox2013_games} \label{prop:pc_passivity}
  The EDM \eqref{eq:pc_dyn} is $\delta$-passive and has a strict storage function given by $$\mathcal S (p, x) = \sum_{i=1}^n \sum_{j=1}^n x_i \int_{0}^{p_j-p_i} \varrho_j(s) \,\mathrm ds$$
\end{prop}

\subsubsection{Perturbed Best Response (PBR) Dynamics \cite{hofbauer2007_jet}} \label{sec:pbr}
\begin{align}
  \dot x(t) = C(p(t)) - x(t)  \label{eq:pbr}
\end{align}
The function $C:\mathbb R^n \to \mathbb X$ is defined as $C(p) = \operatorname*{arg\,max}_{y \in \mathrm{int}(\mathbb X)} \left( p^Ty - v(y) \right)$ where $v:\mathbb X \to \mathbb R$ is a $\mathcal C^2$ function satisfying the following two conditions:\footnote{Note that there is a unique $y$ in $\mathrm{int}(\mathbb X)$ for which $p^Ty - v(y)$ is maximized.}
\begin{align*}
  &z^T \nabla_x^2 v(x) z > 0 ~\text{for all $(x,z)$ in $\mathbb X \times \mathbb{TX} \setminus \{\mathbf 0\}$} \\
  &\left\| \nabla_x v(x) \right\| \to \infty ~ \text{as $x \to \mathrm{bd}(\mathbb X)$}
\end{align*}
We refer to such $C$ and $v$ as the \textit{choice function} and \textit{(deterministic) perturbation}, respectively. Note that \eqref{eq:logit} is a particular case of \eqref{eq:pbr} with the perturbation ${v (x) = \eta \sum_{i=1}^n x_i \ln x_i}$. The following proposition establishes $\delta$-passivity of \eqref{eq:pbr}.
\begin{prop} \label{prop:pbr_passivity}
  The EDM \eqref{eq:pbr} is $\delta$-passive and has a strict storage function given by
  \begin{align} \label{eq:storage_function_pbr}
    \mathcal S(p, x) = \max_{y \in \mathrm{int}(\mathbb X)} \left( p^Ty - v(y)\right) - \left( p^Tx - v(x)\right)
  \end{align}
  If $v$ satisfies the following strong convexity condition then the EDM is strictly output $\delta$-passive with index $\eta > 0$:
  \begin{align*}
    z^T \nabla_x^2 v(x) z \geq \eta z^T z ~ \text{for all $(x,z)$ in $\mathbb X \times \mathbb{TX}$}
  \end{align*}
\end{prop}

The proofs of Propositions \ref{prop:rep_dyn}, \ref{prop:ept_passivity}, and \ref{prop:pbr_passivity} are given in Appendix~\ref{proof:passivity_edm}, whereas the proof of Proposition \ref{prop:pc_passivity} is presented in \cite{fox2013_games}.

\section{Properties of $\delta$-Passive Evolutionary Dynamic Models} \label{sec:properties}
\subsection{Payoff Monotonicity} \label{sec:payoff_monotonicity}
Let us consider the following two conditions \cite{hofbauer2009_jet, hofbauer2011_te} -- \textit{Nash Stationarity} \textbf{(NS)} and \textit{Positive Correlation} \textbf{(PC)}:
\begin{align}
  & \mathcal V \left( p,x \right) = \mathbf 0 \text{ if and only if $x \in \argmax_{y \in \mathbb X} p^Ty$} \tag{\textbf{NS}} \label{eq:ns} \\
  & p^T \mathcal V(p, x) \geq 0 \text{ holds for all } (p, x) \in \mathbb{R}^{n} \times \mathbb X \tag{\textbf{PC}} \label{eq:pc}
\end{align}

We refer to EDMs satisfying both \eqref{eq:ns} and \eqref{eq:pc} as \textit{payoff monotonic}. Examples of payoff monotonic EDMs are those for the EPT dynamics \eqref{eq:ept_dyn} and pairwise comparison dynamics \eqref{eq:pc_dyn}. To provide an intuition behind payoff monotinicity, consider a payoff monotonic EDM with a constant payoff $p(t) = p_0$ for all $t$ in $\mathbb R_+$. The population state $x(t)$ induced by the EDM evolves along the direction $\dot x(t)$ that monotonically increases the average payoff, and it becomes stationary when it reaches the state that attains the maximum average payoff. 

The following two propositions characterize important properties of $\delta$-passive EDMs in connection with payoff monotonicity.
\begin{prop} \label{prop:storage_function_global_minimum}
  For a given $\delta$-passive EDM, let $\mathcal S$ be its storage function and $\mathbb S = \{ (p,x) \,|\, \mathcal V (p,x) = 0 \}$ be its equilibrium points. It holds that $\{ (p,x) \,|\, \mathcal S (p,x) = 0\} \subseteq \mathbb S$ where the equality holds if the EDM satisfies \eqref{eq:ns}.
\end{prop}

\begin{prop} \label{prop:sop_ns_pc}
  Suppose that $n \geq 3$. No EDM can be both strictly output $\delta$-passive and payoff monotonic.
\end{prop}

The proofs of Propositions \ref{prop:storage_function_global_minimum} and \ref{prop:sop_ns_pc} are given in Appendix \ref{proof:properties_passive_dynamics}.

The following corollary is a direct consequence of Proposition \ref{prop:sop_ns_pc}.
\begin{cor} \label{cor:sop_ns_pc}
  The EDMs \eqref{eq:ept_dyn} and \eqref{eq:pc_dyn} are not strictly output $\delta$-passive.
\end{cor}

The proof of Corollary \ref{cor:sop_ns_pc} directly follows from Proposition \ref{prop:sop_ns_pc} and the fact that \eqref{eq:ept_dyn} and \eqref{eq:pc_dyn} are payoff monotonic.

In Section \ref{sec:numerical_ex}, we demonstrate asymptotic stability of equilibrium points of $\delta$-passive EDMs in population games where we will observe that  strictly output $\delta$-passive EDMs achieve the stability in a larger class of population games than do (ordinary) $\delta$-passive EDMs. Besides, the payoff monotonicity ensures that the population state $x(t)$ evolves toward Nash equilibria. Based on these observations, strict output $\delta$-passivity and payoff monotonicty are desired attributes of EDMs; however, Proposition \ref{prop:sop_ns_pc} states that these two attributes are not compatible.

\subsection{Effect of Total Payoff Perturbation} \label{sec:payoff_perturbations}
Motivated by our analysis on the PBR dynamics in Section~\ref{sec:pbr}, we investigate the effect of perturbation on $\delta$-passivity of EDMs. For this purpose, consider the \textit{total payoff function} $u: \mathbb R^n \times \mathbb X \to \mathbb R$ defined by
\begin{align} \label{eq:total_payoff}
  u(p, x) = p^Tx - v(x)
\end{align}
where $v:\mathbb X \to \mathbb R$ is a $\mathcal C^2$ function satisfying the following two conditions:
\begin{subequations} \label{eq:perturbation_cond}
  \begin{align}
    &z^T \nabla_x^2 v(x) z \geq \eta z^T z ~ \text{for all $(x,z)$ in $\mathbb X \times \mathbb{TX}$} \\
    &\left\| \nabla_x v(x) \right\| \to \infty ~ \text{as $x \to \mathrm{bd}(\mathbb X)$}
  \end{align}
\end{subequations}
where $\eta$ is a positive constant. We refer to $v$ as \textit{deterministic perturbation} \cite{hofbauer2002_econometrica} or \textit{control cost} \cite{mattsson2002_geb}. Notice that without the perturbation ($v = 0$), the total payoff \eqref{eq:total_payoff} coincides with the average payoff $p^Tx$.\footnote{The idea of imposing perturbations on the average payoff appeared in game theory and economics to investigate the effect of random perturbations or disutility on choice models \cite{hofbauer2002_econometrica, mattsson2002_geb, hofbauer2007_jet}, to model human choice behavior \cite{mcfadden1974_frontiers_econometrics}, and to analyze the effect of social norms in economic problems \cite{lindbeck1999_qje}.}

In what follows, we investigate the effect of perturbation using a certain class of EDMs whose state-space representation is given as follows:\footnote{We note that the EDM of the PBR dynamics is not of the form \eqref{eq:pedm}, and our analysis given in Section \ref{sec:payoff_perturbations} complements the $\delta$-passivity analysis on the PBR dynamics.}
\begin{align} \label{eq:pedm}
  \dot x_i(t) &= \sum_{j=1}^n x_j(t) \widetilde{\varrho}_{ji}(\nabla_x^T u(p(t),x(t)) \xi_{ji}(x(t))) \nonumber \\
                &\quad - x_i(t) \sum_{j=1}^n \widetilde{\varrho}_{ij}(\nabla_x^T u(p(t),x(t)) \xi_{ij}(x(t)))
\end{align}
where the function $\xi_{ji}:\mathbb X \to \mathbb{TX}$ assigns a tangent vector to each $x$ in $\mathbb X$,  and $\widetilde{\varrho}_{ji}: \mathbb R^n \to \mathbb R^n$ is called a revision protocol which depends on a directional derivative of the total payoff. We refer to \eqref{eq:pedm} as \textit{perturbed} if the perturbation $v$ satisfies the conditions \eqref{eq:perturbation_cond}, and as \textit{unperturbed} if $v = 0$.

The following are examples of \eqref{eq:pedm}.
\begin{align} \label{eq:pbnn}
  \dot x_i(t) &= \left[ \nabla_x^T u(p(t), x(t)) \left( e_i - x(t) \right) \right]_+ \nonumber \\
              & \quad - x_i(t) \sum_{j=1}^n \left[ \nabla_x^T u(p(t), x(t)) \left( e_j - x(t) \right) \right]_+
\end{align}
\begin{align} \label{eq:psmith}
  \dot x_i(t) &= \sum_{j=1}^n x_j(t) \left[ \nabla_x^T u(p(t), x(t)) \left( e_i - e_j \right) \right]_+ \nonumber \\
              & \quad - x_i(t) \sum_{j=1}^n \left[ \nabla_x^T u(p(t), x(t)) \left( e_j - e_i \right) \right]_+
\end{align}
Note that the unperturbed EDMs \eqref{eq:pbnn} and \eqref{eq:psmith} coincide with \eqref{eq:bnn} and \eqref{eq:smith}, respectively. Thus, according to Propositions \ref{prop:ept_passivity} and \ref{prop:pc_passivity}, with $v=0$, \eqref{eq:pbnn} and \eqref{eq:psmith} are $\delta$-passive EDMs. For the case where $v$ satisfies \eqref{eq:perturbation_cond}, the perturbed EDMs \eqref{eq:pbnn} and \eqref{eq:psmith} are strictly output $\delta$-passive with index $\eta > 0$.

In the following proposition, we generalize this observation using \eqref{eq:pedm}.
\begin{prop} \label{prop:payoff_perturbation}
  Given an EDM \eqref{eq:pedm}, suppose that the EDM is $\delta$-passive when $v=0$. If $v$ satisfies \eqref{eq:perturbation_cond}, then the perturbed EDM \eqref{eq:pedm} is strictly output $\delta$-passive with index $\eta > 0$.
\end{prop}

The proof is given in Appendix \ref{proof:prop_payoff_perturbation}.

\section{Stability of $\delta$-passive Evolutionary Dynamic Models} \label{sec:numerical_ex}
In this section, we investigate how $\delta$-passivity is connected to asymptotic stability in population games.

\subsection{Stability in Population Games} \label{sec:passivity_equivalence}
Consider the following closed-loop configuration of a static payoff function \eqref{eq:static_payoff} and $\delta$-passive EDM \eqref{eq:edm}:
\begin{align} \label{eq:closed_loop_01}
  \dot x(t) = \mathcal V(\mathcal F(x(t)), x(t))
\end{align}
where the $\mathcal C^1$ mapping $\mathcal F: \mathbb X \to \mathbb R^n$ satisfies $z^TD\mathcal F(x) z \leq \epsilon z^Tz$ with $\epsilon \geq 0$ for all $(x,z)$ in $\mathbb X \times \mathbb {TX}$, and the storage function $\mathcal S$ of the EDM is strict. The following proposition establishes asymptotic stability of the equilibrium points $\{x \in \mathbb X \,|\, \mathcal V (\mathcal F(x), x) = \mathbf 0\}$ of \eqref{eq:closed_loop_01}.

\begin{prop} \label{prop:stability_static_function}
  Consider the closed-loop \eqref{eq:closed_loop_01}. The equilibrium points of \eqref{eq:closed_loop_01} are asymptotically stable if either of the following conditions holds:
  \begin{enumerate}
  \item The EDM \eqref{eq:edm} is $\delta$-passive and $\epsilon = 0$.

  \item The EDM \eqref{eq:edm} is strictly output $\delta$-passive with index $\eta$ and it holds that $\eta > \epsilon > 0$.
  \end{enumerate}
\end{prop}

The proof of Proposition \ref{prop:stability_static_function} is given in Appendix \ref{proof:passivity_stability_equivalence}. Note that under \eqref{eq:ns}, Proposition \ref{prop:stability_static_function} essentially establishes asymptotic stability of Nash equilibria of $\mathcal F$.

Proposition \ref{prop:stability_static_function} shows that $\delta$-passivity is a sufficient condition for stability in population games with static payoff functions \eqref{eq:static_payoff}. In what follows, we investigate whether $\delta$-passivity is also a necessary condition for stability. To this end, let us consider the following closed-loop configuration of a variant of the smoothed payoff \eqref{eq:smoothed_payoff} and an EDM \eqref{eq:edm} satisfying \eqref{eq:ns}:
\begin{subequations} \label{eq:closed_loop_02}
  \begin{align}
    \dot p(t) &= \widetilde{\mathcal F}(x(t)) - \frac{1}{n} \mathbf 1 \mathbf 1^T p(t) \label{eq:closed_loop_02_a}\\
    \dot x(t) &= \mathcal V(p(t), x(t)) \label{eq:closed_loop_02_b}
  \end{align}
\end{subequations}
where $\widetilde{\mathcal F}:\mathrm{int}(\mathbb X) \to \mathbb R^n$ is defined by $\widetilde{\mathcal F}(x) = \nabla_x \left( f(x) + \sum_{i=1}^n x_i \ln\frac{a_i}{x_i}\right)$ with a $\mathcal C^1$ concave function ${f:\mathbb X \to \mathbb R}$ and positive constants $\{a_i\}_{i=1}^n$.

Note that $\widetilde{\mathcal F}$ can be interpreted as a perturbed payoff function in a concave potential game in which the perturbation is described by $\sum_{i=1}^n x_i \ln\frac{a_i}{x_i}$. At a stationary point $(p,x)$ of \eqref{eq:closed_loop_02_a}, the population state $x$ is a Nash equilibrium of the game\footnote{Note that due to the perturbation $\sum_{i=1}^n x_i \ln\frac{a_i}{x_i}$, Nash equilibrium is in the interior of $\mathbb X$.} which satisfies the following relation: $x \in \mathbb X$ is a stationary point of \eqref{eq:closed_loop_02_a} if and only if
\begin{align*}
  x_i>0 \text{ implies } i \in \argmax_{j \in \{1, \cdots, n\}} \widetilde{\mathcal F}_j(x)
\end{align*}

In the following proposition, we show that $\delta$-passivity of EDMs is equivalent to Lyapunov stability in population games with the smoothed payoff \eqref{eq:closed_loop_02_a}.
\begin{prop} \label{prop:passivity_closed-loop_stability}
  Consider the closed-loop \eqref{eq:closed_loop_02}. The EDM \eqref{eq:closed_loop_02_b} is $\delta$-passive and its storage function is strict if and only if for every smoothed payoff \eqref{eq:closed_loop_02_a}, there exists an energy function $\mathcal E: \mathbb R^n \times \mathbb X \to \mathbb R_+$ for which the following holds:
  \begin{enumerate}
  \item $\mathcal E(p,x) = \mathcal S(p,x) + \max_{y \in \mathrm{int}(\mathrm{\mathbb X})} \widetilde f (y) - \widetilde f(x)$.
    
  \item $\frac{\mathrm d}{\mathrm dt} \mathcal E(p(t), x(t)) \leq 0$ where the equality holds if and only if $x(t) \in \argmax_{y \in \mathbb X} p^T(t)y$.
  \end{enumerate}
  where $\mathcal S$ is a fixed $\mathcal C^1$ function and $\widetilde f$ is defined by $\widetilde f(x) = f(x) + \sum_{i=1}^n x_i \ln\frac{a_i}{x_i}$ for each $x$ in $\mathbb X$.
\end{prop}

The proof of Proposition \ref{prop:passivity_closed-loop_stability} is given in Appendix \ref{proof:passivity_stability_equivalence}.

In fact, $1)$ and $2)$ of Proposition \ref{prop:passivity_closed-loop_stability} yield asymptotic stability of the Nash equilibrium of $\widetilde{\mathcal F}$. To make the argument simple, let's suppose that $x(t)$ stays in a compact subset of $\mathrm{int}(\mathbb X)$ for all $t$ in $\mathbb R_+$, which ensures that the payoff $p(t)$ also stays in a compact subset of $\mathbb R^n$.

Note that $\mathcal E$ is a Lyapunov function which ensures that the closed-loop \eqref{eq:closed_loop_02} satisfies the Lyapunov stability criterion \cite{khalil2001_prentice_hall}. In particular, if $\mathcal E(p(t),x(t)) = 0$ then $x(t)$ is the Nash equilibrium of $\widetilde{\mathcal F}$ and $p(t) = \mathbf 1 a(t)$ for a real number $a(t)$, and due to the payoff dynamics \eqref{eq:closed_loop_02_a}, $p(t)$ converges to $\widetilde{\mathcal F}(x(t))$. In addition, due to $2)$ of Proposition \ref{prop:passivity_closed-loop_stability}, the value of $\mathcal E(p(t), x(t))$ decreases unless $p(t) = \mathbf 1 a(t)$ for a real number $a(t)$; on the other hand, unless $x(t)$ is the Nash equilibrium of $\widetilde{\mathcal F}$, the payoff dynamics \eqref{eq:closed_loop_02_a} would not allow $p(t) = \mathbf 1 a(t)$ and $\mathcal E(p(t), x(t))$ keeps decreasing. The above two observations yield asymptotic stability of $(\widetilde{\mathcal F}(x^{\mathrm{NE}}), x^{\mathrm{NE}})$ under \eqref{eq:closed_loop_02}  where $x^{\mathrm{NE}}$ is the Nash equilibrium of $\widetilde{\mathcal F}$.

\subsection{Numerical Simulations}
Through numerical simulations with $n=3$, we evaluate stability of EDMs in population games. We consider two scenarios where in the first scenario, we consider the BNN dynamics \eqref{eq:bnn} and the logit dynamics \eqref{eq:logit} in the Hypnodisk game \cite{hofbauer2011_te} whose static payoff function is given by
\begin{align} \label{eq:hypnodisk}
  p(t) &= \mathcal H(x(t)) 
\end{align}
The mapping $\mathcal H = (\mathcal H_1, \mathcal H_2, \mathcal H_3)$  is defined by
\begin{align*}
  &\begin{pmatrix} \mathcal H_1 (x_1, x_2, x_3) \\ \mathcal H_2 (x_1, x_2, x_3) \\ \mathcal H_3 (x_1, x_2, x_3) \end{pmatrix} \nonumber \\
  &= \cos \left( \theta(x_1, x_2, x_3) \right) \begin{pmatrix} x_1 - \frac{1}{3} \\ x_2 - \frac{1}{3} \\ x_3 - \frac{1}{3} \end{pmatrix} \nonumber \\
  &\quad + \frac{\sqrt{3}}{3} \sin \left( \theta(x_1, x_2, x_3) \right) \begin{pmatrix} x_2 - x_3 \\ x_3 - x_1 \\ x_1 - x_2 \end{pmatrix} + \frac{1}{3} \begin{pmatrix} 1 \\ 1 \\ 1 \end{pmatrix} 
\end{align*}
where $\theta(x_1, x_2, x_3) = \pi \left[ 1 - b\left( \sum_{i=1}^3 \left| x_i - \frac{1}{3} \right|^2 \right) \right]$ and $b$ is a function satisfying
\begin{enumerate}
\item $b(r) = 1$ if $r \leq R_I^2$
\item $b(r) = 0$ if $r \geq R_O^2$
\item $b(r)$ is decreasing if $R_I^2 < r < R_O^2$
\end{enumerate}
with $R_I = 0.1$ and $R_O = 0.4$.

\begin{figure*}
  \centering
  \subfigure[Multiple population state trajectories induced by the BNN dynamics with different initial conditions in the Hypnodisk game; and the storage function evaluated along a population state trajectory.]{\includegraphics[scale=0.4]{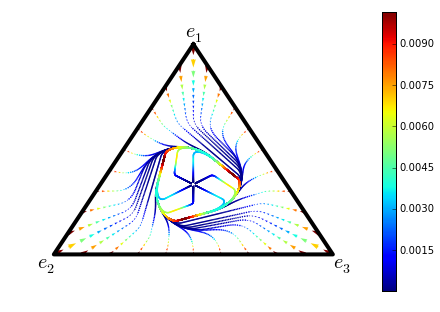} \includegraphics[scale=0.3]{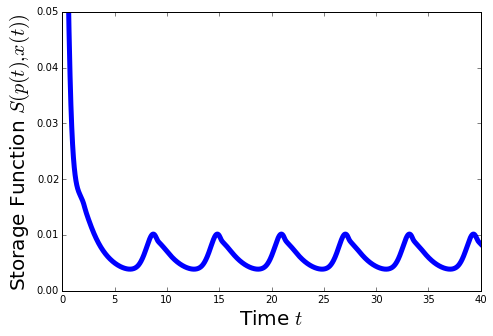}}

  \centering
  \subfigure[Multiple population state trajectories induced by the logit dynamics with different initial conditions in the Hypnodisk game; and the storage function evaluated along a population state trajectory.]{\includegraphics[scale=0.4]{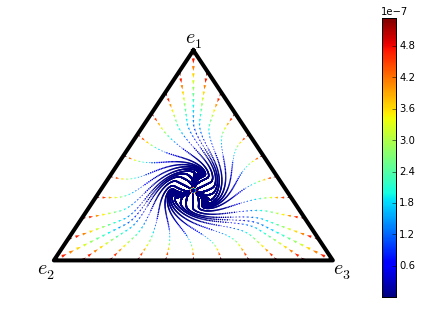} \includegraphics[scale=0.3]{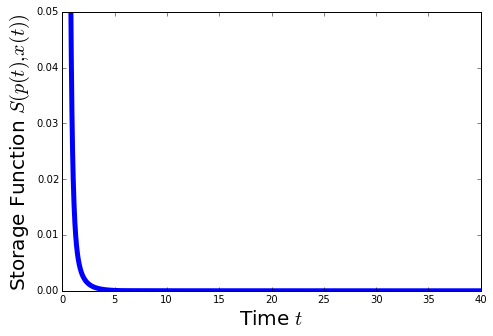}}

  \centering
  \subfigure[A single population state trajectory induced by the replicator dynamics with smoothed payoff; and the KL divergence $\mathcal E(x(t)) = \sum_{i=1}^3 x_i(t) \ln\frac{x_i^{\mathrm{NE}}}{x_i(t)}$ evaluated along the trajectory, where $x^{\mathrm{NE}} = \frac{1}{3}\mathbf 1$.]{\includegraphics[scale=0.4]{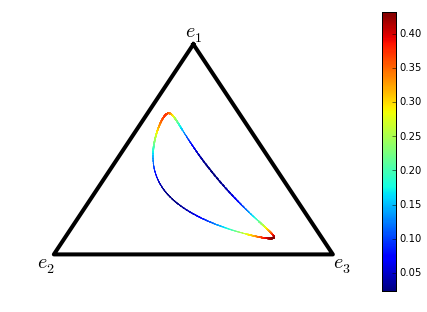} \includegraphics[scale=0.3]{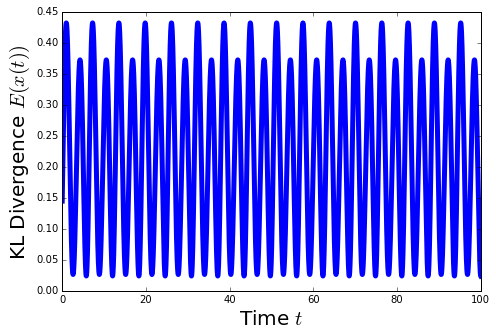}}

  \centering
  \subfigure[A single population state trajectory induced the BNN dynamics with smoothed payoff; and the storage function evaluated along the trajectory.]{\includegraphics[scale=0.4]{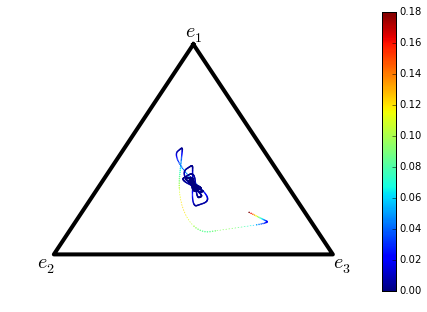} \includegraphics[scale=0.3]{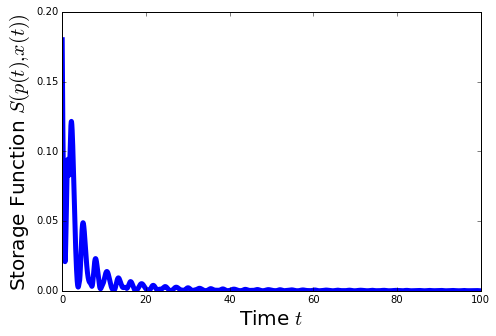}}

  \caption {Simulation results that illustrate asymptotic stability in population games.}
  \label{figure:numerial_simulations}
\end{figure*}

It can be verified that there is a constant $\eta > 0$ for \eqref{eq:logit} for which $2)$ in Proposition \ref{prop:stability_static_function} holds and the population state $x(t)$ converges to  the Nash equilibrium of $$\widetilde{\mathcal H}(x) = \mathcal H(x)+\begin{pmatrix} -\ln x_1 - 1 \\ -\ln x_2 - 1 \\ -\ln x_3 - 1 \end{pmatrix}$$ However, Proposition \ref{prop:stability_static_function} does not apply to \eqref{eq:bnn}. Simulation results given in Fig. \ref{figure:numerial_simulations} (a)-(b) illustrate that the population state trajectories of \eqref{eq:logit} under \eqref{eq:hypnodisk} converge to the Nash equilibrium of $\widetilde{\mathcal H}$ whereas those of \eqref{eq:bnn} do not.

In the second scenario, we compare stability of the replicator dynamics \eqref{eq:rep} and the BNN dynamics \eqref{eq:bnn} using the smoothed payoff \eqref{eq:closed_loop_02_a} where the mapping $\widetilde{\mathcal F}$ is defined by 
$$\widetilde{\mathcal F}(x) = \begin{pmatrix} -x_1 \\ -x_2 \\ -x_3 \end{pmatrix} + \begin{pmatrix} -\ln x_1 - 1 \\ -\ln x_2 - 1 \\ -\ln x_3 - 1 \end{pmatrix}$$
Note that according to Proposition \ref{prop:passivity_closed-loop_stability} and $\delta$-passivity of \eqref{eq:bnn}, we can construct a Lyapunov function $\mathcal E$ and establish asymptotic stability of the Nash equilibrium of $\widetilde{\mathcal F}$ whereas \eqref{eq:rep} is not $\delta$-passive (see Proposition \ref{prop:rep_dyn}) and Proposition~\ref{prop:passivity_closed-loop_stability} does not apply to \eqref{eq:rep}. Simulation results given in Fig. \ref{figure:numerial_simulations} (c)-(d) illustrate that the population state trajectory of \eqref{eq:bnn} under \eqref{eq:closed_loop_02_a} converges to the Nash equilibrium of $\widetilde{\mathcal F}$ whereas that of \eqref{eq:rep} oscillates around the Nash equilibrium.

\section{Conclusions} \label{sec:conclusion}
In this paper, we have investigated $\delta$-passivity of dynamic models in evolutionary game theory. We defined and characterized $\delta$-passivity using the state-space representation of EDMs. Based on the characterization, we studied certain properties of $\delta$-passive EDMs and demonstrated how $\delta$-passivity can be used to establish stability in population games. As future work, it would be interesting to further investigate stability of EDMs in a generalized class of population games in which diverse higher-order dynamics and/or delay are involved.

\begin{appendix}
\renewcommand{\themydef}{\thesection.\arabic{mydef}}
\subsection{Proofs of Theorem \ref{thm:passivity_edm} and Propositions \ref{prop:rep_dyn}, \ref{prop:ept_passivity}, and \ref{prop:pbr_passivity}} \label{proof:passivity_edm}

\subsubsection{Proof of Theorem \ref{thm:passivity_edm}}
  The definition of $\delta$-passivity in Definition \ref{def:passivity} is closely related with the notion of dissipativity in dynamical system theory \cite{willems1972_arma}. To see this, let us rewrite \eqref{eq:edm} in the following form:
  \begin{subequations} \label{eq:edm_expanded}
    \begin{align}
      \dot x(t) &= \mathcal V(u(t), x(t)) \\
      \dot y(t) &= \dot x(t)
    \end{align}
  \end{subequations}
  Note that \eqref{eq:edm_expanded} can be interpreted as a state-space equation for a  control-affine nonlinear system with the input $\dot u(t)$, state $(u(t),x(t))$, and output $\dot y(t)$. According to the definition of dissipativity \cite{willems1972_arma}, the system \eqref{eq:edm_expanded} is dissipative with respect to the supply rate $s(\dot u, \dot y) = \dot u^T \dot y - \eta \dot y^T \dot y$ with a constant $\eta$ if there is a $\mathcal C^1$ function $\mathcal S: \mathbb R^n \times \mathbb X \to \mathbb R_+$ for which
  \begin{align} \label{eq:dissipativity_inequality}
    &\int_{t_0}^{t} \left[ \dot u^T(\tau) \dot y(\tau) - \eta \dot y^T(\tau) \dot y(\tau) \right] \,\mathrm d\tau \nonumber \\
    &\geq \mathcal S(p(t), x(t)) - \mathcal S(p(t_0), x(t_0))
  \end{align}
  holds for all $t \geq t_0 \geq 0$ and $u(\cdot) \in \mathfrak P$. Then, by the equivalence between \eqref{eq:edm} and \eqref{eq:edm_expanded}, we can verify that the $\delta$-passivity condition is satisfied for \eqref{eq:edm} if there is a $\mathcal C^1$ function $\mathcal S$ for which \eqref{eq:dissipativity_inequality} holds with $\eta \geq 0$ for every $u(\cdot)$ in $\mathfrak P$.

Based on the above observation, using the dissipativity characterization theorem (see, for instance, Theorem 1 in \cite{hill1976_ieee_tac}), we can see that there is a $\mathcal C^1$ function $\mathcal S$ satisfying \eqref{eq:p_condition_01} and \eqref{eq:p_condition_02} with $\eta \geq 0$ if and only if under the same choice of $\mathcal S$, the $\delta$-passivity inequality \eqref{eq:passivity_inequality} holds with the same choice of $\eta$ for all $t \geq t_0 \geq 0$ and all $p(\cdot)$ in $\mathfrak P$. \hfill \QED

\subsubsection{Proof of Proposition \ref{prop:rep_dyn}} \label{proof:prop_rep_dyn}
We proceed by showing that any $\mathcal C^1$ function $\mathcal S:\mathbb R^n \times \mathbb X \to \mathbb R_+$ satisfying \eqref{eq:p_condition_01} does not satisfy the condition \eqref{eq:p_condition_02} under \eqref{eq:rep_dyn}. Using Theorem \ref{thm:passivity_edm}, we conclude that \eqref{eq:rep_dyn} is not passive.

Let us re-write \eqref{eq:rep_dyn} in the following form:
\begin{align}
  \dot x_i(t) = \sum_{j=1}^n x_i(t) x_j(t) \left(p_i(t) - p_j(t) \right)
\end{align}

Note that any $\mathcal C^1$ function $\mathcal S$ satisfying \eqref{eq:p_condition_01} should be of the following form:
$$\mathcal S(p, x) = \frac{1}{4}\sum_{i=1}^n \sum_{j=1}^n x_i x_j (p_i - p_j)^2 + g(x)$$
where $g: \mathbb X \to \mathbb R_+$ is a $\mathcal C^1$ function. By taking a partial derivative of $\mathcal S$ with respect to $x$, we obtain \eqref{eq:prop:rep_dyn_01}.
\begin{figure*}
  \begin{align} \label{eq:prop:rep_dyn_01}
    \nabla_x^T \mathcal S (p, x) \mathcal V(p, x) &= \left[ \begin{pmatrix} \frac{1}{2} \sum_{j=1}^n x_j (p_1 - p_j)^2 \\ \vdots \\ \frac{1}{2} \sum_{j=1}^n x_j (p_n - p_j)^2 \end{pmatrix} + \nabla_x g(x) \right]^T \begin{pmatrix} \sum_{j=1}^n x_1 x_j (p_1 - p_j) \\ \vdots \\ \sum_{j=1}^n x_n x_j (p_n - p_j) \end{pmatrix}
  \end{align}
\end{figure*}
Let us choose $x_j = 0$ for all $j \geq 3$. Then, we obtain \eqref{eq:prop:rep_dyn_02}.
\begin{figure*}
  \begin{align} \label{eq:prop:rep_dyn_02}
    \nabla_x^T \mathcal S (p, x) \mathcal V(p, x) &= \begin{pmatrix} \frac{1}{2} x_2 (p_1 - p_2)^2 + \frac{\partial g}{\partial x_1} (x) \\ \frac{1}{2}x_1 (p_1 - p_2)^2 + \frac{\partial g}{\partial x_2} (x) \end{pmatrix} ^T \begin{pmatrix} x_1 x_2 (p_1 - p_2) \\ -x_1 x_2 (p_1 - p_2) \end{pmatrix} \nonumber \\
                                                  &= -\frac{1}{2} x_1 x_2 (p_1 - p_2) \left[ (x_1 - x_2) (p_1 - p_2)^2 + 2 \left( \frac{\partial g}{\partial x_2} (x) - \frac{\partial g}{\partial x_1} (x)\right) \right]
  \end{align}

  \hrulefill
\end{figure*}
Note that for fixed $x_1, x_2$ (except for the points at which $\nabla_x^T \mathcal S(p, x) \mathcal V(p, x) = 0$ holds for all $p$ in $\mathbb R^n$), there exists $p \in \mathbb R^n$ for which $\nabla_x^T \mathcal S (p, x) \mathcal V(p, x) > 0$ holds. Therefore, the function $\mathcal S$ does not satisfy \eqref{eq:p_condition_02}. Since we made an arbitrary choice of $\mathcal S$, we conclude that there does not exist a $\mathcal C^1$ function $\mathcal S$ that satisfies \eqref{eq:p_condition_01} and \eqref{eq:p_condition_02} simultaneously. \hfill \QED

\subsubsection{Proof of Proposition \ref{prop:ept_passivity}} \label{proof:prop_ept_passivity}
We first note that the condition \eqref{eq:acuteness} implies the so-called \textit{Strict Positive Correlation} \eqref{eq:strict_positive_correlation} \cite{sandholm2005_jet} described as follows:
\begin{align}
  \mathcal V(p, x) \neq \mathbf 0 \text{ implies } p^T \mathcal V(p, x) > 0 \tag{\textbf{SPC}} \label{eq:strict_positive_correlation}
\end{align}

Let $\gamma:\mathbb R^n \to \mathbb R$ be a $\mathcal C^1$ function for which \eqref{eq:integrability} holds. It can be verified that $\gamma$ satisfies
\begin{subequations}
  \begin{align}
    \nabla_p \gamma(\hat p) &= \mathcal V(p, x) \label{eq:prop_ept_passivity_04} \\
    \nabla_x^T \gamma(\hat p) \mathcal V(p, x) &= -\left( \mathbf 1^T \varrho(\hat p)\right) \left( p^T \mathcal V(p, x)\right) \label{eq:prop_ept_passivity_01}
  \end{align}
\end{subequations}
where $\hat p = p - p^Tx \mathbf 1$. Let us select a candidate storage function as $\mathcal S(p, x) = \gamma(\hat p)$. Due to \eqref{eq:prop_ept_passivity_04}, the function $\mathcal S$ satisfies \eqref{eq:p_condition_01}. In conjunction with the fact that $\varrho(\hat p) = \mathbf 0$ implies $\mathcal V(p, x) = \mathbf 0$, due to \eqref{eq:strict_positive_correlation} and \eqref{eq:prop_ept_passivity_01}, we can see that \eqref{eq:p_condition_02} holds with $\eta = 0$ and the equality in \eqref{eq:p_condition_02} holds only if $\mathcal V(p,x) = \mathbf 0$.

Suppose that $\gamma$ also satisfies the following inequality for every $\hat p$ in $\mathbb R^n$:
\begin{align} \label{eq:prop_ept_passivity_02}
  \gamma(\hat p) \geq \gamma(\mathbf 0)
\end{align}
Then, without loss of generality by setting $\gamma(\mathbf 0) = 0$, we conclude that $\mathcal S(p,x) = \gamma(\hat p)$ is non-negative and in conjunction with aforementioned arguments, the EDM \eqref{eq:ept_dyn} is $\delta$-passive where its storage function is strict and is given by $\mathcal S(p,x) = \gamma(\hat p)$. In what follows, we show that \eqref{eq:prop_ept_passivity_02} is valid for every $\hat p$ in $\mathbb R^n$.

We first claim that \eqref{eq:prop_ept_passivity_02} holds for all $(p,x)$ in the set $\mathbb S$ of equilibrium points of \eqref{eq:ept_dyn}. To see this, first note that due to \eqref{eq:strict_positive_correlation}, it holds that $\mathcal V(\mathbf 0, x) = \mathbf 0$, i.e., $\mathbf 0 \in \mathbb S_x$, for all $x$ in $\mathbb X$. By \eqref{eq:prop_ept_passivity_04}, for fixed $x$ in $\mathbb X$, the following equality holds for all $p$ in $\mathbb R^n$:

\begin{align}
  \gamma \left( \hat p\right) = \gamma \left( \mathbf 0 \right) + \int_0^1 \dot{\mathbf p}^T(s) \mathcal V(\mathbf p(s), x) \,\mathrm d s
\end{align}
where $\mathbf p: [0,1] \to \mathbb R^n$ is a parameterization of a piece-wise smooth path from $\mathbf{0}$ to $p$ and $\hat p = p- p^Tx \mathbf 1$. According to the path-connectedness assumption (Assumption \ref{assump:path_connected}), for each $p$ in $\mathbb S_x$, there is a path from $\mathbf 0$ to $p$ in which the entire path is contained in $\mathbb S_x$, i.e., $\mathcal V(\mathbf p(s), x) = \mathbf 0$ holds for all $s$ in $[0,1]$; hence the following equality holds for every $p$ in $\mathbb S_x$:
\begin{align} \label{eq:prop_ept_passivity_03}
  \gamma \left( \hat p\right) - \gamma \left( \mathbf 0 \right) = \int_0^1 \dot{\mathbf p}^T(s) \mathcal V(\mathbf p(s), x) \,\mathrm ds = \mathbf 0
\end{align}
Since \eqref{eq:prop_ept_passivity_03} holds for every $(p,x)$ in $\mathbb S$, this proves the claim.

To see that \eqref{eq:prop_ept_passivity_02} extends to the entire domain $\mathbb R^n \times \mathbb X$, by contradiction, let us assume that there is $(p', x') \notin \mathbb S$ for which $\mathcal S(p', x') = \gamma(\hat p') < \gamma(\mathbf 0)$ holds. Let $x(t)$ be the population state induced by \eqref{eq:ept_dyn} with the initial condition $x(0) = x'$ and the constant payoff $p(t) = p'$ for all $t$ in $\mathbb R_+$. By \eqref{eq:strict_positive_correlation} and \eqref{eq:prop_ept_passivity_01}, the value of $\mathcal S(p', x(t))$ is strictly decreasing unless $\mathcal V(p', x(t)) = \mathbf 0$. By the hypothesis that $\mathcal S (p', x') < \gamma(\mathbf 0)$ and by \eqref{eq:prop_ept_passivity_03}, for every $(p,x)$ in $\mathbb S$, it holds that $\mathcal S(p',x') < \mathcal S(p,x)$ and the state $\left( p(t), x(t) \right)$ never converges to $\mathbb{S}$. On the other hand, by LaSalle's Theorem \cite{khalil2001_prentice_hall}, since $p(t)$ is constant and the population state $x(t)$ is contained in a compact set, $x(t)$ converges to an invariant subset of $\left\{ x \in \mathbb X \,\big|\, \nabla_x^T \mathcal S(p', x) \mathcal V(p', x) = 0\right\}$. By \eqref{eq:strict_positive_correlation} and \eqref{eq:prop_ept_passivity_01}, the invariant subset is contained in $\mathbb S$. This contradicts the fact that the state $(p(t), x(t))$ does not converge to $\mathbb S$; hence $\gamma(\hat p) \geq \gamma(\mathbf 0)$ holds for all $(p,x)$ in $\mathbb R^n \times \mathbb X$. \hfill \QED

\subsubsection{Proof of Proposition \ref{prop:pbr_passivity}} \label{proof:prop_pbr_passivity}
The analysis used in Theorem 2.1 of \cite{hofbauer2002_econometrica} suggests that the following hold:
\begin{subequations} \label{eq:prop_perturbed_dynamics_sop_01}
  \begin{align}
    z^T \nabla_p \left[ \max_{y \in \mathrm{int}(\mathbb{X})} (p^Ty - v(y)) \right] = z^T C(p)
  \end{align}
  \begin{align}
    z^T \nabla_x v(x) = z^T p \text{ if and only if } x = C(p) \label{eq:prop_perturbed_dynamics_sop_01_b}
  \end{align}
\end{subequations}
for all $(p,x,z)$ in $\mathbb R^n \times \mathbb X \times \mathbb{TX}$. Using \eqref{eq:prop_perturbed_dynamics_sop_01}, we can see that
\begin{align}
  \nabla_p \mathcal S(p,x) = C(p)-x = \mathcal V(p,x)
\end{align}
and
\begin{align} \label{eq:prop_perturbed_dynamics_sop_02}
  \nabla_x^T \mathcal S(p, x) \mathcal V(p, x) &= - \left( p - \nabla_x v(x) \right)^T \mathcal V(p, x) \nonumber \\
                                      &= -\left( \nabla_y v(y) - \nabla_x v(x)\right)^T \left( y - x \right)
\end{align}
where $y = C(p)$. By the fact that $v$ is strictly convex, it holds that ${\nabla_x^T \mathcal S(p, x) \mathcal V(p, x) \leq 0}$ where the equality holds only if $\mathcal V(p,x) = \mathbf{0}$. According to Theorem \ref{thm:passivity_edm}, we conclude that the EDM \eqref{eq:pbr} is $\delta$-passive and its storage function is strict and is given by \eqref{eq:storage_function_pbr}.

Furthermore, if the perturbation $v$ satisfies $z^T \nabla_x^2 v(x) z \geq \eta z^Tz$ for all $(x,z)$ in $\mathbb X \times \mathbb{TX}$, i.e., $v$ is strongly convex, then it holds that $\left( \nabla_y v(y) - \nabla_x v(x)\right)^T \left( y - x \right) \geq \eta \left\| y-x\right\|^2$ for all $x,y$ in $\mathbb X$. Using \eqref{eq:prop_perturbed_dynamics_sop_02} we can derive $\nabla_x^T \mathcal S(p, x) \mathcal V(p, x) \leq -\eta \mathcal V^T(p, x) \mathcal V(p, x)$. Hence, by Theorem~\ref{thm:passivity_edm}, we conclude that \eqref{eq:pbr} is strictly output $\delta$-passive with index $\eta$. \hfill \QED

\subsection{Proofs of Propositions \ref{prop:storage_function_global_minimum} and \ref{prop:sop_ns_pc}} \label{proof:properties_passive_dynamics}
\subsubsection{Proof of Proposition \ref{prop:storage_function_global_minimum}} \label{proof:prop_storage_function_global_minimum}
The first part of the statement directly follows from the condition \eqref{eq:p_condition_01} and the fact that at a global minimizer $(p^\ast,x^\ast)$ of $\mathcal S$, it holds that $\nabla_p \mathcal S(p^\ast,x^\ast) = \mathbf 0$. Now suppose that the EDM satisfies \eqref{eq:ns}. To prove the second statement, it is sufficient to show that at each equilibrium point $(p,x)$ of \eqref{eq:edm}, it holds that $\mathcal S(p,x) = 0$. To this end, let us consider an \textit{anti-coordination game} whose payoff function is given by $\mathcal F_{x_{o}} (x) = -(x - x_{o})$ for a fixed $x_o$ in $\mathbb{X}$. Notice that $x_o$ is a unique Nash equilibrium of the game. In what follows, we show that $\mathcal S(p_o, x_o) = 0$ holds for any choice of $x_o$ from $\mathbb X$ and $p_o$ from $\mathbb S_{x_o}$.

Let $(p^\ast, x^\ast)$ be a global minimizer of $\mathcal S$, i.e., $\mathcal S\left( p^\ast, x^\ast \right) = 0$. By the first part of the statement and \eqref{eq:ns}, we have that $\mathcal V(p^\ast s, x^\ast) = \mathbf 0$ for all $s$ in $\mathbb R_+$; hence it holds that
$$\mathcal S (\mathbf 0, x^\ast) = \mathcal S (p^\ast, x^\ast) - \int_0^1 \left( p^\ast \right)^T  \mathcal V(p^\ast s, x^\ast) \,\mathrm d s = 0$$ 
By the continuity of $\mathcal S$, for each $\epsilon > 0$ there exists $\delta > 0$ for which $\mathcal S \left( \delta \mathcal F_{x_o} (x^\ast), x^\ast \right) < \epsilon$ holds.

According to the conditions \eqref{eq:p_condition_01} and \eqref{eq:p_condition_02} with $\eta = 0$, the following relation holds for every positive constant $\delta$:
\begin{align} \label{eq:prop_storage_function_global_minimum_02}
  &\frac{\mathrm d}{\mathrm d t} \mathcal S (\delta \mathcal F_{x_o}(x(t)), x(t)) \nonumber \\
  &\leq \delta \mathcal V^T(\delta \mathcal F_{x_o}(x(t)), x(t)) D\mathcal F_{x_0}(x(t)) \mathcal V(\delta \mathcal F_{x_o}(x(t)), x(t)) \nonumber \\
  &= -\delta \left\| \mathcal V \left( \delta \mathcal F_{x_o}(x(t)), x(t) \right) \right\|^2
\end{align}
Suppose that the population state $x(t)$ induced by \eqref{eq:edm} in the anti-coordination game starts from $x(0) = x^\ast$. By an application of LaSalle's theorem \cite{khalil2001_prentice_hall} and by \eqref{eq:ns}, we can verify that $(\delta \mathcal F_{x_o}(x(t)), x(t))$ converges to $\left( \mathbf{0}, x_{o} \right)$ as $t \to \infty$. In addition, due to \eqref{eq:prop_storage_function_global_minimum_02}, we have that 
$$\mathcal S \left( \mathbf 0, x_o \right) \leq \mathcal S \left( \delta \mathcal F_{x_0} (x^\ast), x^\ast \right) < \epsilon$$
Since this holds for every $\epsilon > 0$, we conclude that $\mathcal S (\mathbf 0, x_o) = 0$. By the fact that $\mathcal V (p_o s, x_o) = \mathbf 0$ for all $s$ in $\mathbb R_+$ if $p_o$ belongs to $\mathbb S_{x_o}$, we can see that the following equality holds for every $p_o$ in $\mathbb S_{x_o}$:
\begin{align} \label{eq:prop_storage_function_global_minimum_03}
  \mathcal S (p_o, x_o ) = \mathcal S(\mathbf 0, x_o ) + \int_0^1 p_o^T \mathcal V (p_o s, x_o) \,\mathrm d s = 0
\end{align}

Since we made an arbitrary choice of $x_o$ from $\mathbb X$ in constructing the anti-coordination game, we conclude that \eqref{eq:prop_storage_function_global_minimum_03} holds for every $(p_o, x_o)$ in $\mathbb S$. This proves the proposition. \hfill \QED

\subsubsection{Proof of Proposition \ref{prop:sop_ns_pc}} \label{proof:prop_sop_ns_pc}
We first construct a game $\mathcal F_\nu$ based on the \textit{Hypnodisk game} \cite{hofbauer2011_te}, which is described by the following payoff function: For $a_1, a_2, a_3$ in $\mathbb R_+$, 
\begin{align*}
  &\begin{pmatrix} \mathcal H_1 (a_1, a_2, a_3) \\ \mathcal H_2 (a_1, a_2, a_3) \\ \mathcal H_3 (a_1, a_2, a_3) \end{pmatrix} \nonumber \\
  &= \cos \left( \theta(a_1, a_2, a_3) \right) \begin{pmatrix} a_1 - \frac{1}{3} \\ a_2 - \frac{1}{3} \\ a_3 - \frac{1}{3} \end{pmatrix} \nonumber \\
  &\quad + \frac{\sqrt{3}}{3} \sin \left( \theta(a_1, a_2, a_3) \right) \begin{pmatrix} a_2 - a_3 \\ a_3 - a_1 \\ a_1 - a_2 \end{pmatrix} + \frac{1}{3} \begin{pmatrix} 1 \\ 1 \\ 1 \end{pmatrix} 
\end{align*}
where $\theta(a_1, a_2, a_3) = \pi \left[ 1 - b\left( \sum_{i=1}^3 \left| a_i - \frac{1}{3} \right|^2 \right) \right]$ and $b$ is a bump function that satisfies
\begin{enumerate}
\item $b(r) = 1$ if $r \leq R_I^2$
\item $b(r) = 0$ if $r \geq R_O^2$
\item $b(r)$ is decreasing if $R_I^2 < r < R_O^2$
\end{enumerate}
with $0<R_I<R_O<\frac{1}{\sqrt{6}}$. 

Now consider a payoff function $\mathcal F' = \begin{pmatrix} \mathcal F_1' & \cdots & \mathcal F_n' \end{pmatrix}^T$ defined as follows: For each $i$ in $\{1, \cdots, n\}$,
\begin{align} \label{eq:hypnodisk_extended}
  \mathcal F_i'(x) = \begin{cases} 
    \mathcal H_j (x_1, x_2, x_3 + \cdots + x_n) & \text{ if } j \in \{1, 2\} \\
    \mathcal H_3 (x_1, x_2, x_3 + \cdots + x_n) & \text{ otherwise}
  \end{cases}
\end{align}
Note that the set of Nash equilibria of $\mathcal F'$ is given as follows:
\begin{align} \label{eq:ne_hypnodisk_extended}
  \left\{ x \in \mathbb X \,\bigg|\, x_1 = x_2 = x_3+\cdots+x_n = \frac{1}{3} \right\}
\end{align}
Since $\theta$ is a smooth function, $\mathcal F'$ is continuously differentiable and its differential map $D \mathcal F'$ is bounded, i.e., for some $\delta > 0$, it holds that $z^T D \mathcal F'(x) z < \delta z^Tz$ for all $(x,z)$ in $\mathbb X \times \mathbb{TX}$. Finally, for a given constant $\nu>0$, we define a new payoff function by $\mathcal F_\nu = \frac{\nu}{\delta} \mathcal F'$.

Using the payoff function $\mathcal F_\nu$, we prove the statement of the proposition. By contradiction, suppose that there is an EDM that is both strictly output $\delta$-passive and payoff monotonic. By definition, the EDM satisfies the $\delta$-passivity condition with index $\eta > 0$. Under the payoff function $\mathcal F_\nu$ for which $\nu < \eta$ holds, the time-derivative of the storage function $\mathcal S( \mathcal F_\nu(x(t)), x(t))$ satisfies the following inequality:
\begin{align*}
  \frac{\mathrm d}{\mathrm d t} \mathcal S (\mathcal F_\nu (x(t)), x(t)) \leq - (\eta - \nu) \left\| \mathcal V (\mathcal F_\nu(x(t)), x(t)) \right\|^2
\end{align*}
By an application of LaSalle's theorem \cite{khalil2001_prentice_hall} and by \eqref{eq:ns}, we can verify that the population state trajectory $x(\cdot)$ induced by the EDM under $\mathcal F_\nu$ converges to the set of Nash equilibria \eqref{eq:ne_hypnodisk_extended} as $t \to \infty$.

On the other hand, when $x(t)$ is contained in the set 
\begin{align*}
  &\left\{ x \in \mathbb X \,\Bigg|\, \left( x_1 - \frac{1}{3} \right)^2 + \left( x_2 - \frac{1}{3} \right)^2 \right. \nonumber \\
  & \qquad\qquad \left. + \left( x_3 + \cdots + x_n - \frac{1}{3} \right)^2 \leq R_I^2 \right\}
\end{align*}
by \eqref{eq:pc}, it holds that
\begin{align}
  & \mathcal F_\nu^T(x(t)) \mathcal V ( \mathcal F_\nu(x(t)), x(t)) \nonumber \\
  & = \frac{\nu}{2\delta} \frac{\mathrm d}{\mathrm d t} \left[ \left( x_1(t) - \frac{1}{3} \right)^2 \right. \nonumber \\
  &\quad + \left.\left( x_2(t) - \frac{1}{3} \right)^2 + \left( x_3(t) + \cdots + x_n(t) - \frac{1}{3} \right)^2 \right] \geq 0
\end{align}
Hence, the population state $x(t)$ never converges to the set of Nash equilibria. This is a contradiction; hence, EDMs cannot be both strictly output $\delta$-passive and payoff monotonic. \hfill \QED

\subsection{Proof of Proposition \ref{prop:payoff_perturbation}} \label{proof:prop_payoff_perturbation}
Since the EDM \eqref{eq:pedm} is $\delta$-passive if $v=0$, by Theorem~\ref{thm:passivity_edm}, we can find a storage function $\mathcal S: \mathbb R^n \times \mathbb X \to \mathbb R_+$ for which the conditions \eqref{eq:p_condition_01} and \eqref{eq:p_condition_02} hold for $\eta = 0$. In what follows, we show that when the perturbation $v$ satisfies \eqref{eq:perturbation_cond}, the resulting perturbed EDM is strictly output $\delta$-passive with a storage function $\widetilde{\mathcal S} (p, x) \overset{\mathrm{def}}{=} \mathcal S \left( p - \nabla_x v(x), x \right)$.

By definition, we have that
\begin{align}
  \nabla_{\breve p_i} \mathcal S \left( \breve p,  x \right) &= \sum_{j=1}^n x_j \widetilde{\varrho}_{ji}(\breve p^T \xi_{ji}(x))  - x_i \sum_{j=1}^n \widetilde{\varrho}_{ij}(\breve p^T \xi_{ij}(x)) \nonumber \\
                                                             &\overset{\mathrm{def}}{=} \mathcal V_i(\breve p,x)
\end{align}
Note that \eqref{eq:pedm} can be concisely written as
$$\dot x_i(t) = \mathcal V_i(\breve p,x) \,\big|_{\breve p = p - \nabla_x v(x)}$$

Using $\mathcal V(\breve p,x) = (\mathcal V_1(\breve p,x) \, \cdots \, \mathcal V_n(\breve p,x))^T$, we can compute the gradient of $\widetilde{\mathcal S}$ with respect to $p$ and $x$ as follows:
\begin{align}
  \nabla_p \widetilde{\mathcal S} \left( p, x \right) &= \nabla_{\breve p} \mathcal S \left( \breve p, x \right)\Big|_{\breve p = p - \nabla_x v(x)} \nonumber \\
                                                      &= \mathcal V\left(p - \nabla_x v(x), x \right) \\
  \nabla_x \widetilde{\mathcal S}(p, x) &= \nabla_x^T \left(p - \nabla_x v(x)\right) ~ \nabla_{\breve p} \mathcal S (\breve p, x)\Big|_{\breve p = p-\nabla_x v(x)} \nonumber \\
                                                      &\quad + \nabla_x \mathcal S (\breve p, x) \Big|_{\breve p = p-\nabla_x v(x)} \nonumber \\
                                                      &= - \left( \nabla_x^2 v(x) \right)^T \mathcal V\left(p - \nabla_x v(x), x \right) \nonumber \\
                                                      &\quad + \nabla_x \mathcal S(\breve p, x) \Big|_{\breve p = p-\nabla_x v(x)} \label{eq:S_tilde}
\end{align}
Using \eqref{eq:S_tilde}, we can derive the following:
\begin{align} \label{eq:p_condition_perturbed}
  &\nabla_x^T \widetilde{\mathcal S} (p, x) \mathcal V\left(p - \nabla_x v(x), x \right) \nonumber \\
  &= - \mathcal V^T\left(p - \nabla_x v(x), x \right) \nabla_x^2 v(x) \mathcal V\left(p - \nabla_x v(x), x \right) \nonumber \\
  &\quad + \nabla_x^T \mathcal S (\breve p, x) \mathcal V(\breve p, x) \Big|_{\breve p = p-\nabla_x v(x)} \nonumber \\
  &\leq - \mathcal V^T\left(p - \nabla_x v(x), x \right) \nabla_x^2 v(x) \mathcal V\left(p - \nabla_x v(x), x \right)
\end{align}
where the inequality holds due to $\delta$-passivity of the unperturbed EDM. Since $v$ satisfies $z^T \nabla_x^2 v(x) z \geq \eta  z^Tz$ for all $(x,z)$ in $\mathbb X \times \mathbb{TX}$ with $\eta > 0$, from \eqref{eq:p_condition_perturbed}, we can see that the following inequality holds:
\begin{align*}
  &\nabla_x^T \widetilde{\mathcal S}(p, x) \mathcal V\left(p - \nabla_x v(x), x \right) \\
  &\leq -\eta \mathcal V^T\left(p - \nabla_x v(x), x \right) \mathcal V \left(p - \nabla_x v(x), x \right)
\end{align*}
Hence, using Theorem~\ref{thm:passivity_edm}, we conclude that the perturbed EDM is strictly output $\delta$-passive and its storage function is given by $\widetilde{\mathcal S}$.\hfill \QED

\subsection{Proofs of Propositions \ref{prop:stability_static_function} and \ref{prop:passivity_closed-loop_stability}} \label{proof:passivity_stability_equivalence}

\subsubsection{Proof of Proposition \ref{prop:stability_static_function}}
By taking a time-derivative of the strict storage function $\mathcal S$, we can derive the following relation:
\begin{align}
  &\frac{\mathrm d }{\mathrm d t} \mathcal S(\mathcal F(x(t)), x(t)) \nonumber \\
  &\leq \dot x^T(t) D \mathcal F(x(t)) \dot x(t) \nonumber \\
  &\quad - \eta \mathcal V^T(\mathcal F(x(t)), x(t)) \mathcal V(\mathcal F(x(t)), x(t)) 
\end{align}
Note that according to LaSalle's Theorem \cite{khalil2001_prentice_hall}, if either $1)$ or $2)$ of the statement holds, then the population state $x(t)$ of \eqref{eq:closed_loop_01} converges to the equilibrium points of \eqref{eq:closed_loop_01}. Beside Lyapunov stability directly follows from the fact that under $1)$ or $2)$ of the statement, $\frac{\mathrm d }{\mathrm d t} \mathcal S(\mathcal F(x(t)), x(t)) \leq 0$ holds. \hfill \QED

\subsubsection{Proof of Proposition \ref{prop:passivity_closed-loop_stability}} \label{proof:prop_passivity_closed-loop_stability}
The sufficiency directly follows using Theorem \ref{thm:passivity_edm}, Proposition \ref{prop:storage_function_global_minimum}, and a strict storage function $\mathcal S$ of the EDM. To prove the necessity, we proceed with computing the derivative of the function $\mathcal E$ as follows:

{\small
\begin{align}
  &\frac{\mathrm d}{\mathrm dt} \mathcal E (p(t), x(t)) \nonumber \\
  &= \nabla_p^T \mathcal S(p(t), x(t)) \dot p(t) + \nabla_x^T \mathcal S(p(t), x(t)) \dot x(t) - \widetilde{\mathcal F}^T(x(t)) \dot x(t) \nonumber \\
  &= \nabla_p^T \mathcal S(p(t), x(t)) \left[ \widetilde{\mathcal F}(x(t)) - \frac{1}{n} \mathbf 1 \mathbf 1^T p(t) \right] \nonumber \\
  & \quad + \nabla_x^T \mathcal S(p(t), x(t)) \mathcal V(p(t), x(t)) \nonumber \\
  & \quad - \left[ \widetilde{\mathcal F}(x(t)) - \frac{1}{n} \mathbf 1 \mathbf 1^T p(t) \right]^T \mathcal V(p(t), x(t)) \nonumber \\
  &= \left[ \nabla_p \mathcal S(p(t), x(t)) - \mathcal V(p(t), x(t)) \right]^T \left[ \widetilde{\mathcal F}(x(t)) - \frac{1}{n} \mathbf 1 \mathbf 1^T p(t) \right] \nonumber \\
  & \quad + \nabla_x^T \mathcal S(p(t), x(t)) \mathcal V(p(t), x(t)) \label{eq:passivity_stability_01}
\end{align}}
where we use the fact that $\mathbf 1^T \mathcal V(p, x) = 0$ for all $(p,x)$ in $\mathbb R^n \times \mathbb X$.

Let us select $\widetilde f(x) = \sum_{i=1}^n x_i \ln \frac{a_i}{x_i}$ and define $\overline p(t) = \frac{1}{n} \mathbf 1^T p(t)$. Then, from \eqref{eq:passivity_stability_01}, we can derive the following relation:
\begin{align}
  & \frac{\mathrm d}{\mathrm dt} \mathcal E (p(t), x(t)) \nonumber \\
  &= \left[ \nabla_p \mathcal S(p(t), x(t)) - \mathcal V(p(t), x(t)) \right]^T \begin{pmatrix} \ln \frac{a_1}{x_1}-1 - \overline p(t) \\ \vdots \\ \ln \frac{a_i}{x_i}-1 - \overline p(t) \\ \vdots \\ \ln \frac{a_n}{x_n}-1 - \overline p(t) \end{pmatrix} \nonumber \\
  & \quad + \nabla_x^T \mathcal S(p(t), x(t)) \mathcal V(p(t), x(t)) \nonumber
\end{align}
Recall that the constant $a_i$ can be any positive real number and the image of $\ln a_i $ is $(-\infty, \infty)$. Based on this observation for any $(p(t), x(t))$ in $\mathbb R^n \times \mathrm{int}(\mathbb X)$, we can select the constants $\{a_i\}_{i=1}^n$ for which $\frac{\mathrm d}{\mathrm dt} \mathcal E (p(t), x(t)) > 0$ unless $\nabla_p \mathcal S(p(t), x(t)) = \mathcal V(p(t), x(t))$. Hence, by the continuity of $\nabla_p \mathcal S$ and $\mathcal V$, we have that $\nabla_p \mathcal S(p,x) = \mathcal V(p,x)$ for all $(p,x)$ in $\mathbb R^n \times \mathbb X$. According to \eqref{eq:ns} and $2)$ of Proposition \ref{prop:passivity_closed-loop_stability}, $\nabla_x^T \mathcal S(p, x) \mathcal V(p, x) = 0$ if and only if $\mathcal V(p,x) = \mathbf 0$. Using Theorem \ref{thm:passivity_edm}, we conclude that the EDM is $\delta$-passive with a strict storage function $\mathcal S$. \hfill \QED

\end{appendix}

\bibliographystyle{IEEEtran}
\bibliography{IEEEabrv,population_games_passivity}

\end{document}